\documentclass{article}
\usepackage{authblk}
\usepackage[T2A]{fontenc}
\usepackage[utf8]{inputenc}
\usepackage{amsmath}
\usepackage{amsthm}
\usepackage{amsfonts}
\usepackage{amssymb}
\usepackage{enumitem}
\usepackage[pdftex]{graphicx, color}
\usepackage{hyperref}
\hypersetup{
    colorlinks=true,
    linkcolor=blue,
    filecolor=magenta,
    urlcolor=blue,
}

\renewcommand{\phi}{\varphi}

\newcommand{\Vol}{\mathrm{Vol}\,}

\newcommand{\Did}{\Delta^{id}}
\newcommand{\Dreg}{\Delta^{reg}}
\newcommand{\Vid}{V^{id}}
\newcommand{\Vreg}{V^{reg}}
\newcommand{\lamid}{\lambda^{id}}
\newcommand{\lamreg}{\lambda^{reg}}

\theoremstyle{plain}
\newtheorem{theorem}{Theorem}[section]
\newtheorem{corollary}[theorem]{Corollary}
\newtheorem{lemma}[theorem]{Lemma}
\newtheorem{proposition}[theorem]{Proposition}
\newtheorem{remark}[theorem]{Remark}

\title{\textbf{Volume asymptotics for a family of hyperbolic 3-manifolds with totally geodesic boundary}}
\author[]{D.\,D.~Nigmedyanov}
\affil[]{Saint Petersburg State University, Universitetskaya nab., 7--9, Saint Petersburg 199034, Russia}
\affil[ ]{\textit{e-mail: danil.nig1@gmail.com}}

\date{November 8, 2026}

\begin{document}

\maketitle

\begin{abstract}
For the family $\{\mathcal{M}_{g,k}\}_{g\geqslant k\geqslant 0}$ of orientable hyperbolic 3-manifolds with totally geodesic boundary and cusps introduced by R.~Frigerio, B.~Martelli, and C.~Petronio, the volume of the manifolds depends only on the non-negative integer parameters $g$ and $k$; thus the volume function $\Vol(g,k)$ is well defined. The main result of the paper is an asymptotic description of this function: it is approximated by a linear function of $g$ and $k$ whose coefficients are explicitly expressed in terms of the volumes of the regular ideal tetrahedron and the regular ideal octahedron; moreover, the volume does not exceed the linear approximation by a difference of order $O(1/g)$. As consequences, three theorems comparing volumes are obtained: monotonicity under coordinatewise increase of the parameters, estimates for the passage from parameter $g$ to $g+1$, and monotonicity along lines of constant complexity $g+k$. For small values of the parameters, approximate values of the volumes are computed and collected in a table used in verifying the cases $g<8$.
\end{abstract}

\tableofcontents

\section{Introduction}

Computing the volume of hyperbolic 3-manifolds goes back to the work of N.~I.~Lobachevsky \cite{Lobachevsky}, who obtained an explicit formula for the volume of an ideal tetrahedron in terms of the function $\Lambda$ he introduced. A topological \emph{ideal triangulation} of a compact connected 3-manifold \(M\) is a decomposition of the complement $M\setminus\partial M$ into tetrahedra with vertices removed. For manifolds admitting an ideal triangulation, the volume is computed as the sum of the volumes of the constituent tetrahedra; the systematic development of this approach is due to W.~Thurston \cite{Thurston}, who also established that the set of volumes of complete finite-volume hyperbolic 3-manifolds is well ordered and has order type $\omega^\omega$.

In the \emph{cusped} case, the SnapPea census, certified by B.~Burton \cite{Burton}, provides a complete list of cusped hyperbolic 3-manifolds of small complexity together with their volumes. Beyond tabulation, one of the few infinite series for which the volume is known explicitly is given by \emph{tetrahedral} manifolds, glued from $k$ regular ideal tetrahedra: for these the volume equals $k v_3$, where $v_3\approx 1.014941$ is the volume of the regular ideal tetrahedron. Such a triangulation is minimal by Thurston's bound $c(M)\geqslant \Vol(M)/v_3$, so for tetrahedral manifolds the volume determines the complexity. A census of all cusped tetrahedral manifolds was compiled by E.~Fominykh, S.~Garoufalidis, M.~Goerner, V.~Tarkaev, and A.~Vesnin \cite{Fominykh}.

For manifolds with totally geodesic boundary, the first classification results were given by M.~Fujii \cite{Fujii}, who described the manifolds decomposing into two regular truncated tetrahedra of dihedral angle $\pi/6$; these manifolds were shown by S.~Kojima and Y.~Miyamoto \cite{KojimaMiyamoto} to have the minimum volume among all compact hyperbolic 3-manifolds with nonempty totally geodesic boundary. For manifolds of small complexity, R.~Frigerio, B.~Martelli, and C.~Petronio \cite{FMP04} enumerated all orientable finite-volume hyperbolic 3-manifolds with totally geodesic boundary admitting an ideal triangulation with at most four tetrahedra and computed their volumes.

All the results mentioned above concern either manifolds of bounded complexity (finite classifications, censuses) or tetrahedral manifolds, whose volume is determined by the number of tetrahedra. In contrast, the present paper studies a family of manifolds of \emph{unbounded} complexity, and the question is how the volume behaves as the complexity grows.

In \cite{FrigMarPet03-2}, R.~Frigerio, B.~Martelli, and C.~Petronio described the family $\{\mathcal{M}_{g,k}\}_{g\geqslant k\geqslant 0}$ of compact orientable 3-manifolds with boundary. By definition, a compact orientable 3-manifold $M$ belongs to $\mathcal{M}_{g,k}$ if it admits an ideal triangulation with $g+k$ tetrahedra and its boundary $\partial M$ consists of a closed orientable surface of genus $g$ and $k$ tori. In \cite{FrigMarPet03-2} it is proved that the manifolds in the family $\{\mathcal{M}_{g,k}\}_{g\geqslant k\geqslant 0}$ are hyperbolic, and their volume depends only on the values of non-negative integer parameters $g$ and $k$, so that the volume function $\Vol(g,k)$ is well defined. It is also proved that the sets $\mathcal{M}_{g,k}$ are nonempty for all $g>k\geqslant 0$ with $g\geqslant 2$, and for $g=k>0$ (the latter for even $k$); in particular, $\mathcal{M}_{1,0}=\varnothing$\footnote{The exclusion of $(g,k)=(1,0)$ is due to the orientability condition: if it were dropped, this set would consist of the Gieseking manifold, whose volume equals $v_3$.}. A pair $(g,k)$ of non-negative integers for which $\mathcal{M}_{g,k}\neq\varnothing$ will be called \emph{admissible}.

The main result of the paper is a two-sided estimate for $\Vol(g,k)$ with an explicit $O(1/g)$ remainder. While the volume of the manifolds in $\mathcal{M}_{g,0}$ is known to be minimal for given genus, as follows from Miyamoto's work \cite{Miyamoto}, the asymptotic behaviour of $\Vol(g,k)$ for $k>0$ has not been described previously; the estimate below covers all admissible pairs with $g\ge 8$ and describes the behaviour of the volume for arbitrarily large complexity. Recall that $v_3\approx 1.014941$ and $v_8\approx 3.663864$ denote the volumes of the regular ideal tetrahedron and the regular ideal octahedron, respectively.

\begin{theorem}
    \label{thm:Structure}
    For every admissible pair $(g,k)$ with $g\geqslant 8$, the estimate
    $$
        5v_3 k + v_8(g-k)-\frac{C_0}{g}
        \;\leqslant\;
        \Vol(g,k)
        \;\leqslant\;
        5v_3 k + v_8(g-k),
    $$
    holds, where $C_0=32\pi^2/135$. In particular, $\Vol(g,k)=5v_3\,k + v_8\,(g-k) + O(1/g)$ as $g\to\infty$, and the sign of the correction is always nonpositive.
\end{theorem}

In addition, approximate values of the volumes have been computed for small values of the parameters; they are used in the proof of the corollaries to handle the cases $g<8$.

\begin{theorem}
    \label{thm:volumes}
    Approximate values of the volumes $\Vol(g,k)$ for all admissible pairs $(g,k)$ with $2 \leqslant g \leqslant 15$ are presented in Table~\ref{tab:volumes}.
\end{theorem}

Three qualitative properties of the volume function follow directly from these theorems. First, the volume increases monotonically under coordinatewise increase of the parameters.

\begin{corollary}[Coordinatewise growth of volume]
    \label{cor:VolumeGrowth}
    If for two admissible pairs $(g_1,k_1)$ and $(g_2,k_2)$ we have simultaneously $g_1\leqslant g_2$ and $k_1\leqslant k_2$, then $\Vol(g_2,k_2)\geqslant\Vol(g_1,k_1)$. Moreover, this inequality is strict if $(g_1,k_1)\neq(g_2,k_2)$.
\end{corollary}

When the admissible pairs $(g,k)$ are ordered lexicographically and we pass from level $g$ to level $g+1$, the following inequalities for the volumes hold.

\begin{corollary}[Step in $g$]
    \label{cor:VolumeStepUp}
    For all $g\geqslant 4$ and all $k$ such that the pairs $(g,k)$, $(g+1,k-3)$, and $(g+1,k-2)$ are admissible,
    $$
    \Vol(g+1,k-3) < \Vol(g,k) < \Vol(g+1,k-2).
    $$
\end{corollary}

Comparison of volumes of manifolds at one level of complexity deserves separate attention. As shown in \cite{FrigMarPet03-2}, the complexity of every manifold in $\mathcal{M}_{g,k}$ equals $g+k$.

\begin{corollary}[Same complexity]
    \label{cor:VolumeSameComplexity}
    If $g_1+k_1=g_2+k_2$ and $k_1<k_2$, then
    $\Vol(g_1,k_1) > \Vol(g_2,k_2)$.
\end{corollary}

\begin{remark}
    The constant $C_0$ is not minimal; it is chosen for convenience in the proofs of the corollaries.
\end{remark}

\section{Outline of the proof}
\label{sec:scheme}

The proof of Theorem~\ref{thm:Structure} proceeds in four steps;
the resulting estimate is then applied to the corollaries.

\paragraph{Step 1. Triangulation and angles.}
A manifold from $\mathcal{M}_{g,k}$ is glued from $g-k$ copies of the regular tetrahedron $\Dreg_\beta$ and $2k$ copies of the tetrahedron $\Did_\alpha$ (Section~\ref{sec:triangulation}), where $(\alpha,\beta)$ is the unique solution of system~\eqref{eq:system2}. As $g\to\infty$, the angles tend to zero at rate $1/g$: $\alpha(g,k)\lesssim 1/g$ and $\beta(g,k)\lesssim 1/g$ (Lemma~\ref{lm:Solutions}). It is this decay that allows one to pass to the limiting case $(\alpha,\beta)=(0,0)$ and estimate the resulting error.

\paragraph{Step 2. Volumes of tetrahedra and monotonicity.}
By Lemma~\ref{lm:Cut}, each tetrahedron decomposes into 6 congruent orthoschemes $T(\alpha,\gamma/2,\pi/3)$, and the volume of an orthoscheme is expressed in terms of the Lobachevsky function (Section~\ref{sec:volumes}). This yields explicit integral formulas for $\Vol(\Did_\alpha)$ and $\Vol(\Dreg_\beta)$ (Lemma~\ref{lm:VolumeFormulas}) and the strict monotonicity of both functions in their respective angle (Proposition~\ref{prop:VolumeProperties}). In particular, the limiting values at zero angle are determined:
$$
    \Vol(\Did_0)=\frac52\,v_3,
    \qquad
    \Vol(\Dreg_0)=v_8.
$$

\paragraph{Step 3. Approximating volume and the remainder.}
Set
$$
    \Vol^{*}(g,k):=2k\,\Vol(\Did_0)+(g-k)\,\Vol(\Dreg_0)
    =v_8(g-k)+5v_3 k,
    \qquad
    \delta\Vol(g,k):=\Vol^{*}(g,k)-\Vol(g,k).
$$
The function $\Vol^{*}$ is linear in $(g,k)$ and therefore easily compared for any pairs. The quantity $\delta\Vol$ is nonnegative by the monotonicity from Step~2 and is the sum of the deviations of the volumes of individual tetrahedra from their limiting values at zero angle.

\paragraph{Step 4. Rate of convergence.}
Each such deviation is quadratic in the angle (Lemma~\ref{lm:VolumeDifferences}), while the angles decay as $1/g$ (Step~1). This yields a remainder estimate uniform in $k$:
$$
    \delta\Vol(g,k)\leqslant\frac{C_0}{g},
    \qquad
    C_0=\frac{32\pi^2}{135}
$$
(Lemmas~\ref{lm:DeltaEstimate}, \ref{lm:UniformEstimate}). Together with $\delta\Vol\geqslant 0$, this gives Theorem~\ref{thm:Structure}: the upper bound follows from the nonnegativity of the remainder, the lower bound from $\delta\Vol\leqslant C_0/g$.

\paragraph{Application to the corollaries.}
For two admissible pairs $(g_1,k_1)$ and $(g_2,k_2)$ with $g_1,g_2\geqslant 8$, the decomposition $\Vol=\Vol^{*}-\delta\Vol$ and the nonnegativity of $\delta\Vol$ imply that if
$$
    \Vol^{*}(g_2,k_2)-\Vol^{*}(g_1,k_1)
    >
    \delta\Vol(g_1,k_1)+\delta\Vol(g_2,k_2),
$$
then $\Vol(g_2,k_2)>\Vol(g_1,k_1)$. The right-hand side is at most $\frac{C_0}{g_1}+\frac{C_0}{g_2}\leqslant \frac{C_0}{4}$, while the left-hand side, by the linearity of $\Vol^{*}$, equals
$$
    v_8\,(g_2-g_1)+(5v_3-v_8)\,(k_2-k_1),
$$
and hence depends only on the differences of the parameters. The verification reduces to comparing two numerical constants: the \emph{linear gap}, depending on $(g_2-g_1,\,k_2-k_1)$, and the \emph{deficit} $C_0/4$. All three corollaries are proved by this scheme (Section~\ref{sec:Consequences}); the boundary cases $g<8$ are verified from Table~\ref{tab:volumes}.

\section{Geometric triangulations}
\label{sec:triangulation}

A \emph{generalized tetrahedron} $\hat{\Delta}$ is a region in three-dimensional hyperbolic space $\mathbb{H}^3$ bounded by four hyperplanes, any two of which intersect in $\mathbb{H}^3$. If a triple of planes intersects in $\mathbb{H}^3$, the intersection point is a vertex of the generalized tetrahedron, which we call \emph{finite}. If the triple of planes intersects on the absolute, we say that $\hat{\Delta}$ has an \emph{ideal} vertex. Otherwise, we say that the triple of planes corresponds to an \emph{ultra-ideal} vertex. For a triple of planes corresponding to an ultra-ideal vertex of $\hat{\Delta}$, one can draw a unique truncation plane perpendicular to each of the three planes (see~\cite[Lemma~2.1]{Kojima}). A \emph{partially truncated tetrahedron} $\Delta$ is a region in $\mathbb{H}^3$ bounded by the four faces of $\hat{\Delta}$ and by the truncation planes at all of its ultra-ideal vertices. We will call \emph{lateral hexagon} and \emph{truncation triangle} the intersection of $\Delta$ with a face of $\hat{\Delta}$ and with the truncation plane to an ultra-ideal vertex of $\hat{\Delta}$. The edges of the truncation triangles, which also belong to the lateral hexagons, will be called \emph{boundary edges}. The other edges of the lateral hexagons will be called \emph{internal edges}. 
As is well known (see \cite{Kojima} for the truncated case and \cite[Theorem~2.2]{FrigPet01} for the general classification), the dihedral angles along the internal edges parametrize a partially truncated tetrahedron.
Moreover, the sum of dihedral angles along internal edges sharing a common ultra-ideal vertex is strictly less than $\pi$; the sum of angles along internal edges sharing a common ideal vertex equals $\pi$; and the sum of angles along internal edges emanating from a finite vertex is greater than $\pi$.

\begin{remark}
\label{rem:ProjectiveModel}
It is often convenient to work with generalized tetrahedra in the projective model of $\mathbb{H}^3$ realized in Minkowski space $M^{3,1}$: a point of $\mathbb{H}^3$ is represented by a vector in $M^{3,1}$, defined up to a positive scalar factor, and a hyperplane by a linear hyperplane in $M^{3,1}$. The open unit ball $B^3\subset M^{3,1}$ models $\mathbb{H}^3$, and the three types of vertices admit a uniform description by the position of the representing vector relative to $\partial B^3$:
\begin{itemize}
    \item a \emph{finite} vertex lies in $B^3$ (the representing vector is timelike);
    \item an \emph{ideal} vertex lies on $\partial B^3$ (the representing vector is lightlike);
    \item an \emph{ultra-ideal} vertex lies outside $\overline{B^3}$ (the representing vector is spacelike).
\end{itemize}
The truncating plane at an ultra-ideal vertex $X$ is then the linear hyperplane $X^\perp$; this uniform description treats all three types of vertices on the same footing and will be used several times below, in particular in the proof of Lemma~\ref{lm:Cut}.
\end{remark}

\begin{figure}[htb]
    \centering
    \includegraphics[scale=0.7]{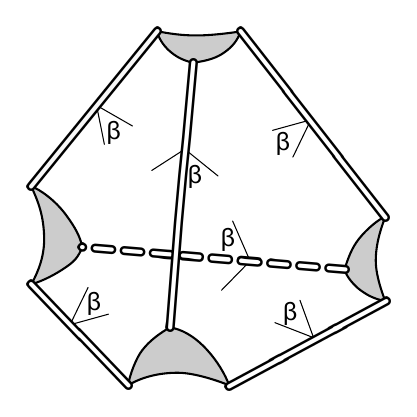}
    \includegraphics[scale=0.7]{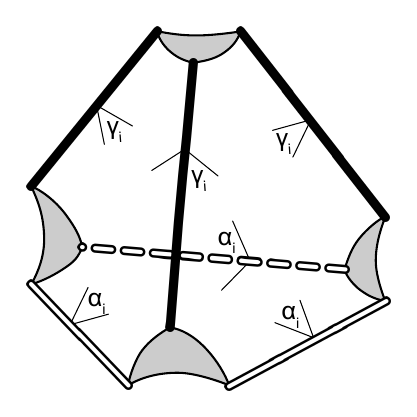}
    \caption{Partially truncated tetrahedra $\Dreg_{\beta}$ (left) and $\Delta_{\alpha,\gamma}$ (right).}
    \label{pic:tetrahedra}
\end{figure}

In~\cite[\S2]{FrigMarPet03-2}, geometric triangulations of the manifolds of the family $\{\mathcal{M}_{g,k}\}_{g\geqslant k \geqslant 0}$ are described. Namely, every manifold in $\mathcal{M}_{g,k}$ can be glued from $g-k$ copies of the regular partially truncated tetrahedron $\Dreg_{\beta}$ with dihedral angle $\beta$ along all internal edges, and from $2k$ copies of the tetrahedron $\Delta_{\alpha,\gamma}$ with dihedral angle $\alpha$ along the internal edges of one lateral hexagon and angle $\gamma$ along the remaining three internal edges (see Fig.~\ref{pic:tetrahedra}).

The hyperbolic structures specified in the interiors of these tetrahedra determine a hyperbolic structure on the entire manifold if and only if the angles of the tetrahedra are solutions of the following system of equations:

\begin{align}
    \label{eq:system}
    \begin{cases}
        6 \gamma = 2\pi \\
        6k\alpha + 6(g-k)\beta = 2 \pi \\
        \frac{\cos^2\beta + \cos\beta}{\sin^2\beta}=\frac{\cos^2\alpha + \cos\gamma}{\sin^2\alpha}.
    \end{cases}
\end{align}

The value $\gamma=\pi/3$ satisfies system~\eqref{eq:system}, so we eliminate $\gamma$ from the system and obtain the following system of equations:

\begin{align}
    \label{eq:system2}
    \begin{cases}
        k\alpha + (g-k)\beta = \pi / 3 \\
        1 - \cos \beta = \frac{2}{3}\sin^2\alpha.
    \end{cases}
\end{align}

Note also that for $\gamma=\pi/3$ the partially truncated tetrahedron $\Delta_{\alpha,\gamma}$ has one ideal vertex. We denote it by $\Delta^{id}_{\alpha}$. According to~\cite{FrigMarPet03-2}, such partially truncated tetrahedra are paired up and, glued in pairs, form the cusps of the hyperbolic manifold under consideration.

\begin{lemma}[Properties of the solution of system~\eqref{eq:system2}]
    \label{lm:Solutions}
    For every admissible pair $(g,k)$ of non-negative integers:
    \begin{enumerate}[label=(\roman*)]
        \item system~\eqref{eq:system2} has a unique solution $(\alpha(g,k),\beta(g,k))\in(0,\pi/3]^2$;
        \item the solution satisfies the inequalities
        $$
        \frac{\sqrt3\,\pi}{6g}\le\alpha(g,k)\le\frac{\pi}{3g},
        \qquad
        \frac{\pi}{3g}\le\beta(g,k)\le\frac{2\pi}{3\sqrt3\,g}.
        $$
    \end{enumerate}
\end{lemma}

\begin{proof}
\emph{Properties of the function $\beta(\alpha)$.}
From the second equation of system~\eqref{eq:system2} we express $\beta$ in terms of $\alpha$:
$$
    1-\cos\beta(\alpha)=\tfrac23\sin^2\alpha,\qquad 0\le \alpha\le\pi/3.
$$
The function $\beta(\alpha)$ is well defined, strictly increasing, concave, with $\beta(0)=0$ and $\beta(\pi/3)=\pi/3$. Differentiating the implicit equation, we obtain
$$
    \sin\beta\cdot\beta'=\tfrac23\sin 2\alpha.
$$
As $\alpha\to0$, we have $\beta\sim\beta'(0)\alpha$, $\sin\beta\sim\beta$, $\sin 2\alpha\sim 2\alpha$. It follows that $\beta'(0)^2=4/3$, i.e., $\beta'(0)=2/\sqrt3$.

From concavity and the boundary values $\beta(0)=0$, $\beta(\pi/3)=\pi/3$, it follows that the graph of $\beta$ lies above the chord:
$$
    \beta(\alpha)\ge\alpha,\qquad 0\le\alpha\le\pi/3.
$$
From concavity and the value $\beta'(0)=2/\sqrt3$, it follows that the graph lies below the tangent at zero:
$$
    \beta(\alpha)\le\frac{2}{\sqrt3}\,\alpha,\qquad 0\le\alpha\le\pi/3.
$$
Thus, for all $0\le\alpha\le\pi/3$,
\begin{equation}
    \label{eq:betaBounds}
    \alpha\le\beta(\alpha)\le\frac{2}{\sqrt3}\,\alpha.
\end{equation}

\emph{Existence and uniqueness.}
Set
$$
    F(\alpha)=k\alpha+(g-k)\beta(\alpha),\qquad 0\le\alpha\le\pi/3.
$$
The function $F$ is continuous and strictly increasing as a linear combination of strictly increasing functions with nonnegative coefficients, and
$$
    F(0)=0,
    \qquad
    F(\pi/3)=g\cdot\frac{\pi}{3}\ge\frac{\pi}{3}.
$$
By the intermediate value theorem, there exists $\bar\alpha\in(0,\pi/3]$ with $F(\bar\alpha)=\pi/3$; by strict monotonicity, it is unique. The pair $(\bar\alpha,\beta(\bar\alpha))$ is the required solution. Part~(i) is proved.

\emph{Estimates.}
Substitute $\alpha=\bar\alpha$, $\beta=\beta(\bar\alpha)$ into the equation $k\alpha+(g-k)\beta=\pi/3$ and use~\eqref{eq:betaBounds}.

From $\beta\ge\alpha$ we obtain $g\bar\alpha \leqslant \pi /3$, so $\bar\alpha\le\pi/(3g)$. From the same inequality $\beta\ge\alpha$ and $k\le g$ we get $\bar\beta\ge\pi/(3g)$.

From $\beta\le(2/\sqrt3)\alpha$ and $k\ge0$:
$$
    \frac{\pi}{3}
    =k\bar\alpha+(g-k)\bar\beta
    \le\Bigl(k+\tfrac{2}{\sqrt3}(g-k)\Bigr)\bar\alpha.
$$
For $0\le k\le g$ the coefficient of $\bar\alpha$ does not exceed $\tfrac{2}{\sqrt3}g$ (the maximum is attained at $k=0$), so
$$
    \frac{\pi}{3}\le\tfrac{2}{\sqrt3}\,g\,\bar\alpha,
    \qquad\text{i.e.,}\qquad
    \bar\alpha\ge\frac{\sqrt3\,\pi}{6g}.
$$
Finally, $\bar\beta\le(2/\sqrt3)\bar\alpha\le 2\pi/(3\sqrt3\,g)$. Part~(ii) is proved.
\end{proof}

\section{Orthoschemes and volumes of tetrahedra}
\label{sec:volumes}

\subsection{Orthoschemes and the Lobachevsky formula}

An \emph{orthoscheme} is a partially truncated tetrahedron whose three internal edges not lying in one face have dihedral angle $\pi/2$. Following~\cite{Kras}, we denote the vertices of the generalized tetrahedron corresponding to the orthoscheme by $v_0,v_1,v_2,v_3$, choosing the order so that the three acute angles are along the edges $v_0v_1$, $v_1v_2$, $v_2v_3$, and the three remaining edges $v_0v_2$, $v_0v_3$, $v_1v_3$ have dihedral angle $\pi/2$. By $T(\phi,\psi,\xi)$ we denote the orthoscheme with $\theta(v_0v_1)=\phi$, $\theta(v_1v_2)=\psi$, $\theta(v_2v_3)=\xi$.

\begin{figure}[htb]
    \centering
    \includegraphics[scale=0.8]{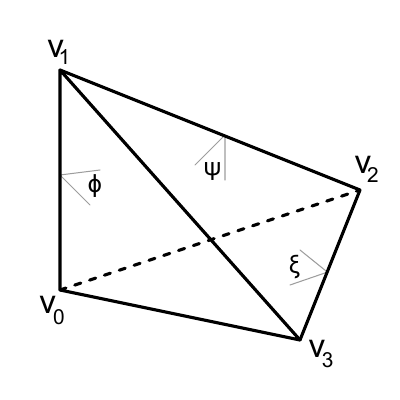}
    \caption{Orthoscheme $T(\phi,\psi,\xi)$.}
    \label{pic:Ortosheme}
\end{figure}

The volume of an orthoscheme is given by the Lobachevsky formula \cite{Lobachevsky}:
\begin{equation}
    \label{eq:VolumeOfOrtoscheme}
    4\Vol(T(\phi,\psi,\xi))
    =
    \Lambda(\phi+\delta)-\Lambda(\phi-\delta)
    +\Lambda(\xi+\delta)-\Lambda(\xi-\delta)
    +\Lambda(\tfrac{\pi}{2}-\psi-\delta)-\Lambda(\tfrac{\pi}{2}-\psi+\delta)
    +2\Lambda(\tfrac{\pi}{2}-\delta),
\end{equation}
where
$$
    \delta=\arctan\!\left(
    \frac{\sqrt{\cos^2\psi-\sin^2\phi\sin^2\xi}}{\cos\phi\cos\xi}
    \right),
$$
and $\Lambda$ is the Lobachevsky function
$$
    \Lambda(x)=\int_0^x -\log|2\sin t|\,dt.
$$

\subsection{Decomposition of tetrahedra into orthoschemes}

It is known (see~\cite[Remark 4.1]{Kras}) that every partially truncated tetrahedron decomposes into 6 orthoschemes. Owing to the symmetries of the tetrahedron $\Delta_{\alpha,\gamma}$, the dihedral angles of all six orthoschemes are expressed in terms of the dihedral angles of the tetrahedron as follows.

\begin{lemma}[Decomposition into orthoschemes]
    \label{lm:Cut}
    The partially truncated tetrahedron $\Delta_{\alpha,\gamma}$ decomposes into 6 congruent orthoschemes $T(\alpha,\gamma/2,\pi/3)$.
\end{lemma}

\begin{proof}
\emph{The case of finite vertices.}
Let $\Delta$ be a hyperbolic tetrahedron with finite vertices $A$, $B$, $C$, $D$ such that the dihedral angles along $AD$, $BD$, $CD$ equal $\gamma$ and those along $AB$, $BC$, $CA$ equal $\alpha$, and let $H$ be the centre of the face $ABC$. 
By the symmetries of $\Delta$, the planes through $DH$ and the medians $AA'$, $BB'$, $CC'$ of $ABC$ cut the faces $BCD$, $CAD$, $ABD$ along the medians from $D$ and are perpendicular to these faces. These three planes divide $\Delta$ into six congruent tetrahedra. By the same symmetries (rotations by $2\pi/3$ about $DH$ and reflections), each of these tetrahedra is an orthoscheme, with three acute dihedral angles $\alpha$, $\gamma/2$, $\pi/3$. 
The verification repeats verbatim the Euclidean construction of the decomposition of a symmetric tetrahedron into orthoschemes: it uses only incidences and perpendicularities and does not rely on Euclid's fifth postulate. Therefore it carries over to $\mathbb{H}^3$ without changes. 
This completes the proof in the case of finite vertices. It remains to treat the general case, where some vertices may be ultra-ideal.

\emph{The general case.}
We now pass to the projective model of $\mathbb{H}^3$ (Remark~\ref{rem:ProjectiveModel}) in order to treat all types of vertices uniformly. In this model a generalized tetrahedron is given by four linearly independent vectors in $M^{3,1}$, and a hyperplane by a linear hyperplane (see~\cite[\S3]{FrigPet01} for an explicit Minkowski implementation of Kojima's truncation, and \cite{Ushijima03} for the projective viewpoint on generalized tetrahedra). All objects of the construction---the points $D,H,A,B,C$ and the midpoints of the sides of $ABC$---are represented by vectors, and the cutting planes are linear hyperplanes in $M^{3,1}$. The construction of the previous paragraph carries over verbatim. 
If a vertex $X$ is ultra-ideal, then the truncating hyperplane $X^\perp$ is perpendicular to every cutting plane passing through $X$; this follows immediately from the definition (the normal to $X^\perp$ is $X$, while the normal to the cutting plane is orthogonal to $X$). Therefore, after truncation of ultra-ideal vertices, the decomposition restricts correctly to the partially truncated tetrahedron $\Delta_{\alpha,\gamma}$ and yields six congruent (partially truncated) orthoschemes. The orthoscheme angles are unchanged by truncation: $\alpha$, $\gamma/2$, $\pi/3$.
\end{proof}

\begin{remark}
\label{rem:IntrinsicCut}
Both the decomposition into orthoschemes and the proof of congruence
of the resulting orthoschemes could be carried out in terms of the
intrinsic geometry of $\mathbb{H}^3$. We choose the Minkowski
description because it treats finite, ideal and ultra-ideal vertices on
an equal footing and avoids a case-by-case discussion of the truncation.
\end{remark}

\begin{corollary}
    \label{cor:CutRegular}
    For $\alpha=\gamma=\beta$, the tetrahedron $\Delta_{\alpha,\gamma}$ coincides with the regular truncated tetrahedron $\Dreg_\beta$ and decomposes into 6 congruent orthoschemes $T(\beta,\beta/2,\pi/3)$.
\end{corollary}

\subsection{Volume formulas}

\begin{lemma}[Volume formulas]
    \label{lm:VolumeFormulas}
    For all $0\leqslant\alpha,\beta\leqslant\pi/3$, the following representations hold:
    \begin{align}
        \Vol(\Did_\alpha)
        &=
        3\Lambda\!\left(\tfrac{\pi}{6}\right)
        +
        \frac32\int_{\alpha-\pi/3}^{\alpha+\pi/3}
        -\log|2\sin t|\,dt,
        \label{eq:VolumeAlpha}\\[2pt]
        \Vol(\Dreg_\beta)
        &=
        8\Lambda\!\left(\tfrac{\pi}{4}\right)
        -
        3\int_0^\beta
        \operatorname{arccosh}\!\left(\frac{\cos t}{2\cos t-1}\right)dt.
        \label{eq:VolumeBeta}
    \end{align}
\end{lemma}

\begin{proof}
    By Lemma~\ref{lm:Cut}, $\Did_\alpha$ decomposes into 6 copies of the orthoscheme $T(\alpha,\pi/6,\pi/3)$. In formula~\eqref{eq:VolumeOfOrtoscheme} for these parameters,
    $$
        \delta
        =
        \arctan\!\left(
        \frac{\sqrt{\cos^2(\pi/6)-\sin^2\alpha\,\sin^2(\pi/3)}}
        {\cos\alpha\,\cos(\pi/3)}
        \right)
        =
        \arctan\sqrt3
        =
        \frac{\pi}{3}.
    $$
    This gives
    \begin{align}
        4\Vol\bigl(T(\alpha,\pi/6,\pi/3)\bigr)
        &=
        \Lambda(\alpha+\pi/3)-\Lambda(\alpha-\pi/3)+2\Lambda(\pi/6) \\
        &=
        2\Lambda\!\left(\tfrac{\pi}{6}\right)
        +
        \int_{\alpha-\pi/3}^{\alpha+\pi/3}
        -\log|2\sin t|\,dt.
    \end{align}
    This yields the volume formula for $\Did_\alpha$. Formula~\eqref{eq:VolumeBeta} is Miyamoto's formula \cite{Miyamoto} for the volume of the regular truncated tetrahedron.
\end{proof}

\subsection{Boundary values and monotonicity}
\label{subsec:VolumeProperties}

Recall the standard notation:
$$
v_3:=3\Lambda\!\left(\tfrac{\pi}{3}\right)\approx 1.014941,
\qquad
v_8:=8\Lambda\!\left(\tfrac{\pi}{4}\right)\approx 3.663864,
$$
the volumes of the regular ideal tetrahedron and the regular ideal octahedron, respectively.

For brevity, when studying volumes as functions of the angle we shall write
\begin{equation}
    \label{eq:VNotation}
    \Vid(\alpha):=\Vol(\Did_\alpha),
    \qquad
    \Vreg(\beta):=\Vol(\Dreg_\beta),
    \qquad
    0\leqslant\alpha,\beta\leqslant\pi/3.
\end{equation}

\begin{proposition}
\label{prop:VolumeProperties}
The functions $\Vid$ and $\Vreg$ on the segment $[0,\pi/3]$ possess the following properties:
\begin{enumerate}
    \item[(i)] boundary values:
    $$
    \Vid(0)=\frac52\,v_3,
    \qquad
    \Vid(\pi/3)=v_3,
    $$
    $$
    \Vreg(0)=v_8,
    \qquad
    \Vreg(\pi/3)=v_3;
    $$

    \item[(ii)] both functions are strictly decreasing, and their derivatives have the form
    \begin{equation}
    \label{eq:VolumeDerivatives}
    \bigl(\Vid\bigr)'(\alpha)
    =
    -3\operatorname{artanh}\!\Bigl(\frac{\tan\alpha}{\sqrt3}\Bigr),
    \qquad
    \bigl(\Vreg\bigr)'(\beta)
    =
    -3\,h(\beta),
    \end{equation}
    where
    $$
    h(\beta):=\operatorname{arccosh}\!\left(\frac{\cos\beta}{2\cos\beta-1}\right).
    $$
\end{enumerate}
\end{proposition}

\begin{proof}
    \emph{Boundary values.}
    Substituting $\alpha=0$ into~\eqref{eq:VolumeAlpha} and using $\Lambda(\pi/6)=\tfrac32\Lambda(\pi/3)$, we obtain
    $$
        \Vid(0)
        =
        3\Lambda(\pi/6) + \frac{3}{2}\cdot 2\Lambda(\pi/3)
        =
        \frac{9}{2}\Lambda(\pi/3) + 3\Lambda(\pi/3)
        =
        \frac{15}{2}\Lambda(\pi/3)
        =
        \frac{5}{2} v_3.
    $$
    For $\beta=0$ the integral in~\eqref{eq:VolumeBeta} vanishes, hence $\Vreg(0) = 8\Lambda(\pi/4) = v_8$. For $\alpha=\beta=\pi/3$ both tetrahedra become regular ideal, so their volumes are $v_3$.

    \emph{Derivatives.}
    Differentiating~\eqref{eq:VolumeAlpha} with respect to $\alpha$ and using $\Lambda'(x) = -\log|2\sin x|$, we obtain
    $$
        \bigl(\Vid\bigr)'(\alpha)
        =
        \frac{3}{2}[\log(2\sin(\pi/3-\alpha)) - \log(2\sin(\pi/3+\alpha))]
        =
        \frac{3}{2}\log\frac{\sin(\pi/3-\alpha)}{\sin(\pi/3+\alpha)}.
    $$
    Substituting $\sin(\pi/3\pm\alpha) = \frac{1}{2}\cos\alpha(\sqrt3\pm\tan\alpha)$, we arrive at
    $$
        \bigl(\Vid\bigr)'(\alpha)
        =
        \frac{3}{2}\log\frac{\sqrt3-\tan\alpha}{\sqrt3+\tan\alpha}
        =
        -3\operatorname{artanh}\left(\frac{\tan\alpha}{\sqrt3}\right).
    $$
    Similarly, from~\eqref{eq:VolumeBeta}:
    $$
        \bigl(\Vreg\bigr)'(\beta)
        =
        -3h(\beta), \qquad h(\beta)
        =
        \operatorname{arccosh}\left(\frac{\cos\beta}{2\cos\beta-1}\right).
    $$
    Both derivatives are negative on $(0,\pi/3)$: the first because $\tan\alpha < \sqrt3$ for $\alpha < \pi/3$, the second because $\cos\beta/(2\cos\beta-1) > 1$ for $\beta < \pi/3$.
\end{proof}

\begin{figure}[htb]
    \centering
    \includegraphics[width=0.6\textwidth]{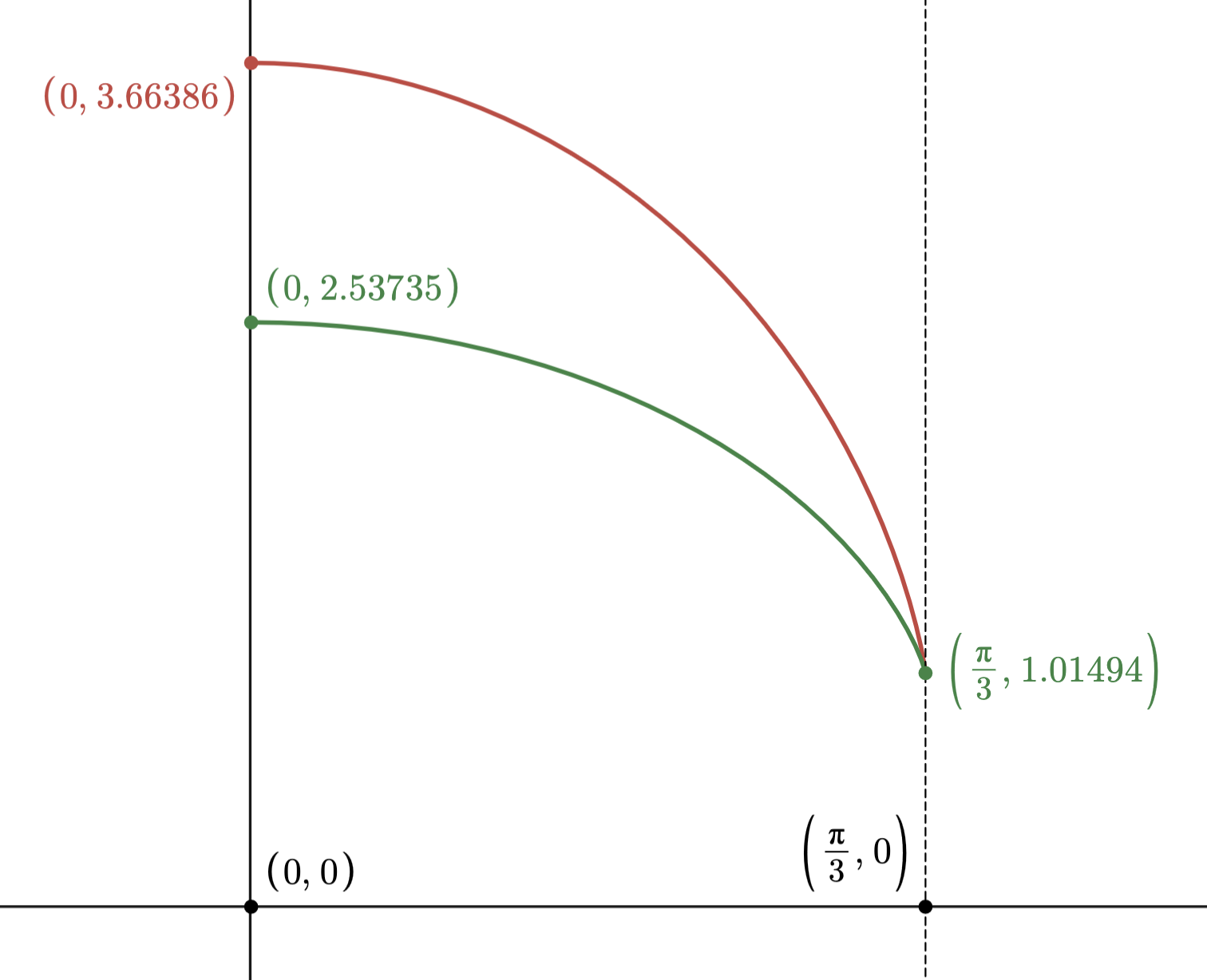}
    \caption{Graphs of $\Vid(\alpha)$ (green) and $\Vreg(\beta)$ (red) on the segment $[0,\pi/3]$.}
    \label{pic:VolumeFunctions}
\end{figure}

\section{Asymptotics of the volume function}
\label{sec:asymptotics}

In this section we approximate the volume $\Vol(g,k)$ by a linear function $\Vol^{*}(g,k)$ and estimate the remainder $\delta\Vol(g,k)$. The main tools are the \emph{deviation functions} $D^{id}, D^{reg}$, which measure the deviation of the volume of an individual tetrahedron from its limiting value at zero angle.

\subsection{Deviation functions}
\label{subsec:DeviationFunctions}

The \emph{deviation functions} are defined by
$$
    D^{id}(\alpha):=\Vid(0)-\Vid(\alpha),
    \qquad
    D^{reg}(\beta):=\Vreg(0)-\Vreg(\beta),
    \qquad
    0\leqslant\alpha,\beta\leqslant\pi/3.
$$

By Proposition~\ref{prop:VolumeProperties}, both functions $\Vid$ and $\Vreg$ are strictly decreasing, so $D^{id}$ and $D^{reg}$ are nonnegative on $[0,\pi/3]$, with $D^{id}(\alpha)=0$ only for $\alpha=0$ and $D^{reg}(\beta)=0$ only for $\beta=0$.

\begin{lemma}
    \label{lm:VolumeDifferences}
    Set
    $$
        \lamid:=\frac{\sqrt3}{2}+0.1,
        \qquad
        \lamreg:=\frac32+0.1.
    $$
    For all \(0\le x\le \pi/6\), the following estimates hold:
    $$
        D^{id}(x)\le \lamid\, x^2,
        \qquad
        D^{reg}(x)\le \lamreg\, x^2.
    $$
\end{lemma}

The constant $0.1$ absorbs the higher-order terms in the Taylor expansions (Appendix~\ref{app:derivatives}); it is chosen for convenience and is not optimal.

\begin{proof}
    Differentiating the definitions of the deviation functions and using formula~\eqref{eq:VolumeDerivatives} (with the opposite sign), we obtain
    \begin{equation}
        \label{eq:DDerivatives}
        \bigl(D^{id}\bigr)'(x)
        =
        -\bigl(\Vid\bigr)'(x)
        =
        3\operatorname{artanh}\!\Bigl(\frac{\tan x}{\sqrt3}\Bigr),
        \qquad
        \bigl(D^{reg}\bigr)'(x)
        =
        -\bigl(\Vreg\bigr)'(x)
        =
        3\,h(x).
    \end{equation}
    Integrating these formulas and expanding the integrands at zero, we arrive at the expansions
    $$
        D^{id}(x)=\frac{\sqrt3}{2}x^2+\frac{x^4}{3\sqrt3}+\frac{2x^6}{15\sqrt3}+R_{id}(x),
        \qquad
        |R_{id}(x)|<1.3\,x^8,
    $$
    $$
        D^{reg}(x)=\frac32x^2+\frac{5}{16}x^4+\frac{25}{192}x^6+R_{reg}(x),
        \qquad
        |R_{reg}(x)|<0.15\,x^8
    $$
    for $0\leqslant x\leqslant\pi/6$; the remainders $R_{id}, R_{reg}$ are estimated by Taylor's formula with the Lagrange form of the remainder in terms of the eighth derivatives on $[0,\pi/6]$ (Appendix~\ref{app:derivatives}). Substituting $x^4\leqslant\pi^2 x^2/36$, $x^6\leqslant\pi^4 x^2/1296$, $x^8\leqslant\pi^6 x^2/46656$ gives
    $$
        D^{id}(x)\leqslant\Bigl(\frac{\sqrt3}{2}+0.1\Bigr)x^2=\lamid\, x^2,
        \qquad
        D^{reg}(x)\leqslant\Bigl(\frac32+0.1\Bigr)x^2=\lamreg\, x^2.
    $$
\end{proof}

\subsection{Limiting volume and estimate of the remainder}

For an admissible pair $(g,k)$, the \emph{limiting volume} is defined by
$$
\Vol^{*}(g,k)
:=
2k\,\Vid(0)+(g-k)\,\Vreg(0)
=
5v_3\,k+v_8\,(g-k).
$$

The function $\Vol^{*}$ is linear in $(g,k)$.  Recall from Section \ref{subsec:DeviationFunctions} that $D^{id}(\alpha)=\Vid(0)-\Vid(\alpha)$ and $D^{reg}(\beta)=\Vreg(0)-\Vreg(\beta)$; hence
\begin{equation}
    \label{eq:VolDecomposition}
    \Vol(g,k)
    =
    \Vol^{*}(g,k)
    -
    2k\,D^{id}\bigl(\alpha(g,k)\bigr)
    -
    (g-k)\,D^{reg}\bigl(\beta(g,k)\bigr).
\end{equation}
By Proposition~\ref{prop:VolumeProperties}, the remainder
$$
    \delta\Vol(g,k):=\Vol^{*}(g,k)-\Vol(g,k)
    =
    2k\,D^{id}\bigl(\alpha(g,k)\bigr)
    +
    (g-k)\,D^{reg}\bigl(\beta(g,k)\bigr)
$$
is nonnegative; this gives the upper bound $\Vol\leqslant\Vol^{*}$ in Theorem~\ref{thm:Structure}. 
Define the majorant of the remainder $\delta\Vol(g,k)$ by the formula
\begin{equation}
    \label{eq:Edef}
    E(g,k):=
    2k\,\lamid\left(\frac{\pi}{3g}\right)^2
    +
    (g-k)\,\lamreg\left(\frac{2\pi}{3\sqrt3\,g}\right)^2.
\end{equation}

\begin{lemma}[Estimate of the remainder]
    \label{lm:DeltaEstimate}
    For all admissible pairs $(g,k)$ with $g\geqslant 8$, the estimate
    $$
        \delta\Vol(g,k)\leqslant E(g,k)
    $$
    holds.
\end{lemma}

\begin{proof}
    By Lemma~\ref{lm:Solutions},
    $$
        \alpha(g,k)\leqslant\frac{\pi}{3g},
        \qquad
        \beta(g,k)\leqslant\frac{2\pi}{3\sqrt3\,g}.
    $$
    For $g\geqslant 8$ this gives
    $$
        \alpha(g,k) \leqslant \frac{\pi}{24} < \frac{\pi}{6},
        \qquad
        \beta(g,k) \leqslant \frac{\pi}{12\sqrt3} < \frac{\pi}{6},
    $$
    so that for $\alpha=\alpha(g,k)$ and $\beta=\beta(g,k)$ Lemma~\ref{lm:VolumeDifferences} is applicable:
    $$
        D^{id}\bigl(\alpha(g,k)\bigr)
        \leqslant \lamid\left(\frac{\pi}{3g}\right)^2,
        \qquad
        D^{reg}\bigl(\beta(g,k)\bigr)
        \leqslant \lamreg\left(\frac{2\pi}{3\sqrt3\,g}\right)^2.
    $$
    Substituting into~\eqref{eq:VolDecomposition}, we obtain $\delta\Vol(g,k)\leqslant E(g,k)$.
\end{proof}

\begin{lemma}[Properties of $E(g,k)$]
    \label{lm:UniformEstimate}
    \begin{enumerate}[label=(\roman*), before=\leavevmode]
        \item The majorant $E(g,k)$ decreases in $k$ for fixed $g$ and decreases in $g$ for fixed $k$ on the domain $0\leqslant k\leqslant g$, $g>0$;
        \item for all $g\geqslant 8$ and $0\leqslant k\leqslant g$,
        $$
            E(g,k)\leqslant \frac{C_0}{g},
            \qquad
            C_0:=\frac{32\pi^2}{135}\approx 2.340.
        $$
    \end{enumerate}
\end{lemma}

\begin{proof}
    We rewrite the formula:
    $$
        E(g,k)
        =
        \frac{\pi^2}{27g^2}\bigl[4\lamreg\,g+k(6\lamid-4\lamreg)\bigr],
        \qquad
        6\lamid-4\lamreg=3\sqrt3-5.8\approx -0.604<0.
    $$

\emph{Monotonicity.}
From the above representation,
$$
	\frac{\partial E}{\partial k}
	=
	\frac{\pi^2(6\lamid-4\lamreg)}{27g^2}<0,
$$
so $E$ decreases in $k$. Furthermore,
$$
	\frac{\partial E}{\partial g}
	=
	\frac{\pi^2}{27g^2}\left[-4\lamreg-\frac{2k(6\lamid-4\lamreg)}{g}\right].
$$
The second summand is positive ($6\lamid-4\lamreg<0$), so we verify
$$
	4\lamreg+\frac{2k(6\lamid-4\lamreg)}{g}
	=
	4\lamreg-\frac{2k|6\lamid-4\lamreg|}{g}
	\geqslant
	4\lamreg-2|6\lamid-4\lamreg|.
$$
Since $2\lamreg=3.2$ and $|6\lamid-4\lamreg|\approx 0.604$, we have
$4\lamreg-2|6\lamid-4\lamreg|\approx 6.4-1.208>0$, whence $\partial E/\partial g<0$.
Part~(i) is proved.

\emph{Uniform estimate.}
From the monotonicity in $k$:
$$
	E(g,k)\leqslant E(g,0)
	=
	\frac{4\lamreg\pi^2}{27g}.
$$
Substituting $\lamreg=8/5$:
$$
	\frac{4\lamreg\pi^2}{27}
	=
	\frac{32\pi^2}{135}
	=
	C_0.
$$
Part~(ii) is proved.
\end{proof}

\subsection{Proof of Theorem~\ref{thm:Structure}}
\label{subsec:ProofStructure}

Let $(g,k)$ be an admissible pair with $g\geqslant 8$.
By the decomposition~\eqref{eq:VolDecomposition},
$$
	\Vol(g,k)=\Vol^{*}(g,k)-\delta\Vol(g,k),
	\qquad
	\Vol^{*}(g,k)=v_8(g-k)+5v_3\,k.
$$

\emph{Upper bound.}
By Proposition~\ref{prop:VolumeProperties}, the functions $\Vreg$ and $\Vid$ are strictly decreasing, so the deviation functions satisfy $D^{id} \geqslant 0$ and $D^{reg} \geqslant 0$. Consequently, $\delta\Vol(g,k)\geqslant 0$, and
$$
	\Vol(g,k)\leqslant \Vol^{*}(g,k)=v_8(g-k)+5v_3\,k.
$$

\emph{Lower bound.}
By Lemma~\ref{lm:DeltaEstimate} and Lemma~\ref{lm:UniformEstimate},
$$
	\delta\Vol(g,k)\leqslant\frac{C_0}{g}.
$$
It follows that
$$
	\Vol(g,k)\geqslant\Vol^{*}(g,k)-\frac{C_0}{g}
	=v_8(g-k)+5v_3\,k-\frac{C_0}{g}.
$$

Both bounds are proved, which completes the proof of Theorem~\ref{thm:Structure}.

\section{Proof of the corollaries}
\label{sec:Consequences}

\subsection{General scheme}
\label{subsec:ComparisonScheme}

Let $(g_1,k_1)$ and $(g_2,k_2)$ be admissible pairs with $g_1,g_2\geqslant 8$, and suppose we wish to show $\Vol(g_2,k_2)>\Vol(g_1,k_1)$.
By the decomposition~\eqref{eq:VolDecomposition} and the nonnegativity of $\delta\Vol$,
$$
	\Vol(g_2,k_2)-\Vol(g_1,k_1)
	\geqslant
	\bigl[\Vol^{*}(g_2,k_2)-\Vol^{*}(g_1,k_1)\bigr]
	-
	\bigl[\delta\Vol(g_1,k_1)+\delta\Vol(g_2,k_2)\bigr].
$$
Therefore, it is enough that
\begin{equation}
	\label{eq:Criterion}
	\Vol^{*}(g_2,k_2)-\Vol^{*}(g_1,k_1)
	>
	\delta\Vol(g_1,k_1)+\delta\Vol(g_2,k_2).
\end{equation}
Using the explicit form $\Vol^{*}(g,k)=v_8(g-k)+5v_3 k$, the left-hand side equals
\begin{equation}
	\label{eq:GapFormula}
	\Vol^{*}(g_2,k_2)-\Vol^{*}(g_1,k_1)
	=
	v_8\,(g_2-g_1)+(5v_3-v_8)\,(k_2-k_1).
\end{equation}
The right-hand side is estimated by Lemma~\ref{lm:UniformEstimate}(ii):
$$
	\delta\Vol(g_1,k_1)+\delta\Vol(g_2,k_2)
	\leqslant
	\frac{C_0}{g_1}+\frac{C_0}{g_2}
	\leqslant
	\frac{2C_0}{8}=\frac{C_0}{4}.
$$
If the quantity obtained from~\eqref{eq:GapFormula} is greater than $C_0/4$, then the inequality $\Vol(g_2,k_2)>\Vol(g_1,k_1)$ is proved.

\subsection{Numerical constants}
A direct computation gives
$$
    v_3\approx 1.0149,\quad
    v_8\approx 3.6639,\quad
    5v_3-v_8\approx 1.4108,
    \qquad
    \frac{C_0}{4}=\frac{8\pi^2}{135}\approx 0.585.
$$
In one case (the left inequality of Corollary~\ref{cor:VolumeStepUp}), the estimate $C_0/4$ turns out to be insufficient, and we use the refined value
$$
    E(8,3)+E(9,0)\approx 0.542
$$
(a direct computation from the definition of the majorant $E(g,k)$).

\subsection{Proof of Corollary~\ref{cor:VolumeGrowth}}
Let $g_2\geqslant g_1$ and $k_2\geqslant k_1$; the case $(g_1,k_1)=(g_2,k_2)$ is trivial (equality). Then at least one of the inequalities $g_1<g_2$ or $k_1<k_2$ is strict, and by~\eqref{eq:GapFormula},
$$
    \Vol^{*}(g_2,k_2)-\Vol^{*}(g_1,k_1)
    =
    v_8(g_2-g_1)+(5v_3-v_8)(k_2-k_1)
    \geqslant
    \min(v_8,5v_3-v_8)
    =
    5v_3-v_8
    \approx
    1.411
    >
    0.585
    =
    \frac{C_0}{4}.
$$
By criterion~\eqref{eq:Criterion}, we obtain $\Vol(g_2,k_2)>\Vol(g_1,k_1)$.

\subsection{Proof of Corollary~\ref{cor:VolumeStepUp}}
\label{subsec:ProofVolumeStepUp}
\emph{Left inequality.}
Set $(g_1,k_1)=(g+1,k-3)$, $(g_2,k_2)=(g,k)$. For $3\leqslant k\leqslant g$ and $g\geqslant 8$, by~\eqref{eq:GapFormula},
$$
    \Vol^{*}(g,k)-\Vol^{*}(g+1,k-3)
    =
    -v_8+3(5v_3-v_8)
    =
    15v_3-4v_8
    \approx
    0.569.
$$
The direct estimate of the deficit $C_0/4\approx 0.585$ is insufficient here, so we refine it using the monotonicity of the majorant $E$ (see definition~\eqref{eq:Edef}) in both arguments (Lemma~\ref{lm:UniformEstimate}(i)):
$$
    \delta\Vol(g+1,k-3)+\delta\Vol(g,k)
    \leqslant
    E(g+1,k-3)+E(g,k)
    \leqslant
    E(9,0)+E(8,3)
    \approx
    0.542
    <
    0.569.
$$
By criterion~\eqref{eq:Criterion}, we obtain $\Vol(g,k)>\Vol(g+1,k-3)$.

\emph{Right inequality.}
Set $(g_1,k_1)=(g,k)$, $(g_2,k_2)=(g+1,k-2)$. For $2\leqslant k\leqslant g$ and $g\geqslant 8$, by~\eqref{eq:GapFormula},
$$
    \Vol^{*}(g+1,k-2)-\Vol^{*}(g,k)
    =
    v_8-2(5v_3-v_8)
    =
    3v_8-10v_3
    \approx
    0.842
    >
    0.585
    =
    \frac{C_0}{4}.
$$
By criterion~\eqref{eq:Criterion}, we obtain $\Vol(g+1,k-2)>\Vol(g,k)$.

\subsection{Proof of Corollary~\ref{cor:VolumeSameComplexity}}
Let $g_1+k_1=g_2+k_2$ and $k_1<k_2$. Then $g_1-g_2=k_2-k_1=:d>0$. By~\eqref{eq:GapFormula},
$$
    \Vol^{*}(g_1,k_1)-\Vol^{*}(g_2,k_2)
    =
    v_8\,d+(5v_3-v_8)(-d)
    =
    d\,(2v_8-5v_3)
    \approx
    2.253\,d
    \geqslant
    2.253
    >
    0.585
    =
    \frac{C_0}{4}.
$$
By criterion~\eqref{eq:Criterion}, we obtain $\Vol(g_1,k_1)>\Vol(g_2,k_2)$.

\subsection{Boundary cases}
In Corollaries~\ref{cor:VolumeGrowth}, \ref{cor:VolumeStepUp}, and~\ref{cor:VolumeSameComplexity}, the main argument assumed $g\geqslant 8$ for both compared pairs. We now consider the cases where this is not so.

\emph{Corollary~\ref{cor:VolumeGrowth}.}
If $g_1<8$ and $g_2<8$, both pairs appear in Table~\ref{tab:volumes}, and the inequality is verified directly. If $g_1<8\leqslant g_2$, we use the intermediate pair $(8,k_1)$, which is admissible since $k_1\leqslant g_1<8$:
$$
    \Vol(g_1,k_1)\overset{\text{table}}{\leqslant}\Vol(8,k_1)
    \leqslant\Vol(g_2,k_2).
$$
The first inequality is read off from the table (both pairs $(g_1,k_1)$ and $(8,k_1)$ have $g\leqslant 8<15$). The second follows from the already proved part for the pairs $(8,k_1)$ and $(g_2,k_2)$ with $g\geqslant 8$, since $g_2\geqslant 8$ and $k_2\geqslant k_1$.

\emph{Corollary~\ref{cor:VolumeStepUp}.}
If $g\geqslant 8$, the main argument works. If $4\leqslant g<8$, all three pairs $(g,k)$, $(g+1,k-3)$, $(g+1,k-2)$ have $g'\leqslant g+1\leqslant 8<15$ and are therefore covered by Table~\ref{tab:volumes}; both inequalities are verified directly.

\emph{Corollary~\ref{cor:VolumeSameComplexity}.}
If $g_1\geqslant 8$, the main argument works. If $g_1<8$, then $g_2=g_1-(k_2-k_1)<g_1<8$, so both pairs are covered by Table~\ref{tab:volumes} and the inequality is verified directly.

\section*{Tool Use}
	
	The numerical computations in Appendix~\ref{app:table} and the verification of the bounds in Appendix~\ref{app:derivatives} were performed using Wolfram Mathematica~14 \cite{Mathematica}. 
	The verification scripts in Appendices~\ref{app:derivatives} and~\ref{app:table} were drafted with the assistance of a large language model and then reviewed, tested, and validated by the authors against
	independent computations. 
	A large language model was also used for assistance with English language editing and \LaTeX{} formatting; all such output was reviewed and verified by the authors, who take full responsibility for the content of the paper.

\appendix

\section{Estimates of higher derivatives}
\label{app:derivatives}

In this section we justify the numerical constants used in the proof of Lemma~\ref{lm:VolumeDifferences}.

\subsection{Setting of the problem}

By formulas~\eqref{eq:VolumeAlpha}--\eqref{eq:VolumeBeta} and Proposition~\ref{prop:VolumeProperties}, the deviation functions have derivatives
$$
    \bigl(D^{id}\bigr)'(x)=3\operatorname{artanh}\!\Bigl(\frac{\tan x}{\sqrt3}\Bigr),
$$
$$
    \bigl(D^{reg}\bigr)'(x) = 3\operatorname{arccosh}\!\left(\frac{\cos x}{2\cos x-1}\right).
$$
Integrating and expanding the integrands at zero, we obtain (see the proof of Lemma~\ref{lm:VolumeDifferences})
\begin{equation}
    \label{eq:AppTaylorId}
    D^{id}(x)=\frac{\sqrt3}{2}x^2+\frac{1}{3\sqrt3}x^4+\frac{2}{15\sqrt3}x^6+R_{id}(x),
\end{equation}
\begin{equation}
    \label{eq:AppTaylorReg}
    D^{reg}(x)=\frac{3}{2}x^2+\frac{5}{16}x^4+\frac{25}{192}x^6+R_{reg}(x).
\end{equation}
The goal is to prove the estimates
\begin{equation}
    \label{eq:AppTarget}
    |R_{id}(x)|<1.3\,x^8,
    \qquad
    |R_{reg}(x)|<0.15\,x^8,
    \qquad
    0\leqslant x\leqslant\pi/6.
\end{equation}

\subsection{Taylor's formula with remainder}

Both functions $D^{id}, D^{reg}$ are analytic on $[0,\pi/6]$: the denominators in the expressions for the derivatives are bounded away from zero,
$$
    4\cos^2 x-1\geqslant 4\cos^2(\pi/6)-1=2,
    \qquad
    2\cos x-1\geqslant 2\cos(\pi/6)-1=\sqrt3-1>0.
$$
Both derivatives are even, so by Taylor's formula with the Lagrange form of the remainder, truncated at the $x^6$ term,
$$
    R_{id}(x) = \frac{(D^{id})^{(8)}(\xi_1)}{8!}x^8,
    \qquad
    R_{reg}(x)=\frac{(D^{reg})^{(8)}(\xi_2)}{8!}x^8
$$
for some $\xi_1,\xi_2\in(0,x)$. It suffices to show that
\begin{equation}
    \label{eq:AppMDefs}
    M_{id}:=\max_{0\leqslant x\leqslant\pi/6}\bigl|(D^{id})^{(8)}(x)\bigr|<5\cdot 10^4,
    \qquad
    M_{reg}:=\max_{0\leqslant x\leqslant\pi/6}\bigl|(D^{reg})^{(8)}(x)\bigr|<6\cdot 10^3.
\end{equation}

\subsection{Estimates of the eighth derivatives}

The maxima in~\eqref{eq:AppMDefs} were found in \texttt{Wolfram Mathematica}~14 using \texttt{NMaximize} with \texttt{WorkingPrecision -> 50}. Both eighth derivatives attain their maximum absolute value at the right endpoint $x=\pi/6$:
$$
    M_{id}=\bigl|(D^{id})^{(8)}(\pi/6)\bigr|\approx 4.1\cdot 10^4
    <5\cdot 10^4,
$$
$$
    M_{reg}=\bigl|(D^{reg})^{(8)}(\pi/6)\bigr|\approx 5.0\cdot 10^3
    <6\cdot 10^3.
$$
The given upper bounds contain a margin of approximately $20\%$ above the certified values. Substituting into~\eqref{eq:AppMDefs} gives
$$
    |R_{id}(x)|\leqslant\frac{5\cdot 10^4}{8!}x^8\approx 1.24\,x^8<1.3\,x^8,
$$
$$
    |R_{reg}(x)|\leqslant\frac{6\cdot 10^3}{8!}x^8\approx 0.149\,x^8<0.15\,x^8,
$$
which matches~\eqref{eq:AppTarget}.

\subsection{Reduction to quadratic estimates}

Using~\eqref{eq:AppTarget} and $x^4\leqslant(\pi^2/36)x^2$, $x^6\leqslant(\pi^4/1296)x^2$, $x^8\leqslant(\pi^6/46656)x^2$ on $[0,\pi/6]$, we obtain
$$
    \frac{x^4}{3\sqrt3}+\frac{2x^6}{15\sqrt3}+|R_{id}(x)|
    \leqslant
    \underbrace{\left(\frac{\pi^2}{108\sqrt3}+\frac{2\pi^4}{19440\sqrt3}+\frac{1.3\,\pi^6}{46656}\right)}_{K_{id}}x^2,
$$
$$
    \frac{5}{16}x^4+\frac{25}{192}x^6+|R_{reg}(x)|
    \leqslant
    \underbrace{\left(\frac{5\pi^2}{576}+\frac{25\pi^4}{248832}+\frac{0.15\,\pi^6}{46656}\right)}_{K_{reg}}x^2.
$$
Numerically,
$$
    K_{id}\approx 0.085<0.1,
    \qquad
    K_{reg}\approx 0.099<0.1.
$$
This suffices to complete the proof of Lemma~\ref{lm:VolumeDifferences}.

\subsection{Verification script}

\begin{verbatim}
gD[x_] := 3 ArcTanh[Tan[x]/Sqrt[3]];
gR[x_] := 3 ArcCosh[Cos[x]/(2 Cos[x] - 1)];

Mid  = NMaximize[{Abs[D[gD[x], {x, 7}]], 0 <= x <= Pi/6}, x,
         WorkingPrecision -> 50];
Mreg = NMaximize[{Abs[D[gR[x], {x, 7}]], 0 <= x <= Pi/6}, x,
         WorkingPrecision -> 50];

Kid  = Pi^2/(108 Sqrt[3]) + 2 Pi^4/(19440 Sqrt[3]) + 1.3 Pi^6/46656;
Kreg = 5 Pi^2/576 + 25 Pi^4/248832 + 0.15 Pi^6/46656;

Print["M_id  = ", Mid[[1]]];
Print["M_reg = ", Mreg[[1]]];
Print["K_id  = ", N[Kid,  20]];
Print["K_reg = ", N[Kreg, 20]];
\end{verbatim}

\section{Table of volumes}
\label{app:table}

\subsection{Method of computation}

For each admissible pair $(g,k)$, we found the solution $(\alpha,\beta)$ numerically by Newton's method, taking as initial approximations the midpoints of the intervals from Lemma~\ref{lm:Solutions}:
$$
    \alpha_0=\frac{\pi}{3g}\cdot\frac{2+\sqrt3}{4},
    \qquad
    \beta_0=\frac{\pi}{3g}\cdot\frac{3+2\sqrt3}{6},
$$
which guarantees convergence of Newton's method to a solution in the admissible domain. We ran the iterations to machine precision and substituted the resulting $\alpha,\beta$ into~\eqref{eq:VolumeAlpha}--\eqref{eq:VolumeBeta}.

No singularities arise in the computation from these formulas. Indeed, for all admissible pairs with $g\geqslant 2$ we have
$$
    \alpha(g,k)\leqslant\frac{\pi}{3g}\leqslant\frac{\pi}{6},
    \qquad
    \beta(g,k)\leqslant\frac{2\pi}{3\sqrt3\,g}\leqslant\frac{\pi}{3\sqrt3}<\frac{\pi}{3},
$$
so the integrand in~\eqref{eq:VolumeBeta} has no singularities on $[0,\beta]$ (its denominator $2\cos t-1$ is bounded away from zero below by $2\cos(\pi/(3\sqrt3))-1>0$). The integral in~\eqref{eq:VolumeAlpha} was computed via the Lobachevsky function using the formula
$$
    \int_{\alpha-\pi/3}^{\alpha+\pi/3}-\log|2\sin t|\,dt
    =
    \Lambda\!\left(\alpha+\tfrac{\pi}{3}\right)
    -
    \Lambda\!\left(\alpha-\tfrac{\pi}{3}\right),
$$
where $\Lambda$ is expressed in terms of the dilogarithm:
$$
    \Lambda(\theta)=\tfrac12\,\mathrm{Cl}_2(2\theta)
    =\tfrac12\,\mathrm{Im}\,\mathrm{Li}_2(e^{2i\theta}).
$$
This avoids numerical integration in~\eqref{eq:VolumeAlpha}.

We performed all computations in \texttt{Wolfram Mathematica}~14 (built-in functions \texttt{FindRoot}, \texttt{NIntegrate}, \texttt{PolyLog}) with 50 significant digits of precision; the table contains rounding to four decimal places.

Table~\ref{tab:volumes} collects the values $\Vol(g,k)$ and the corresponding limiting volumes $\Vol^{*}(g,k)$ for all admissible pairs with $2\leqslant g\leqslant 15$. Direct subtraction confirms the inequality $\Vol(g,k)\leqslant\Vol^{*}(g,k)$ and the decay of the remainder $\Vol^{*}(g,k)-\Vol(g,k)$ as $g$ grows, in accordance with Theorem~\ref{thm:Structure}. The values for $g<8$ are used in Corollaries~\ref{cor:VolumeGrowth}--\ref{cor:VolumeSameComplexity} to verify the boundary cases.

\subsection{Verification script}

\begin{verbatim}
    alpha0[g_] := (2 + Sqrt[3])/12 * Pi/g;
    beta0[g_]  := (3 + 2 Sqrt[3])/18 * Pi/g;

    Lambda[x_] := 1/2 Im[PolyLog[2, Exp[2 I x]]];
    Vid[a_] := 3 Lambda[Pi/6]
        + 3/2 (Lambda[a + Pi/3] - Lambda[a - Pi/3]);
    Vreg[b_] := 8 Lambda[Pi/4] - 3 NIntegrate[
        ArcCosh[Cos[t]/(2 Cos[t] - 1)], {t, 0, b},
        Method -> "LocalAdaptive", WorkingPrecision -> 50];
    solve[g_, k_] := FindRoot[
        {k a + (g - k) b == Pi/3, 1 - Cos[b] == 2/3 Sin[a]^2},
        {{a, alpha0[g]}, {b, beta0[g]}},
        WorkingPrecision -> 50];
    vol[g_, k_] := Module[{s = solve[g, k]},
        2 k Vid[a /. s] + (g - k) Vreg[b /. s]];
\end{verbatim}

\begin{table}[htbp]
\centering
\footnotesize
\setlength{\tabcolsep}{3pt}
\renewcommand{\arraystretch}{0.92}
\begin{tabular}{|cc r r || cc r r || cc r r|}
\hline
$g$ & $k$ & $\Vol(g,k)$ & $\Vol^{*}(g,k)$ & $g$ & $k$ & $\Vol(g,k)$ & $\Vol^{*}(g,k)$ & $g$ & $k$ & $\Vol(g,k)$ & $\Vol^{*}(g,k)$ \\
\hline
2&0&6.4520&7.3277 & 9&3&37.0154&37.2073 & 13&1&48.9130&49.0411 \\
2&1&7.7976&8.7386 & 9&4&38.4232&38.6181 & 13&2&50.3225&50.4519 \\
2&2&9.1345&10.1494 & 9&5&39.8309&40.0290 & 13&3&51.7320&51.8628 \\
3&0&10.4286&10.9916 & 9&6&41.2385&41.4398 & 13&4&53.1414&53.2736 \\
3&1&11.8127&12.4024 & 9&7&42.6460&42.8507 & 13&5&54.5508&54.6844 \\
3&2&13.1942&13.8133 & 9&8&44.0534&44.2615 & 13&6&55.9602&56.0953 \\
4&0&14.2382&14.6555 & 10&0&36.4738&36.6386 & 13&7&57.3696&57.5061 \\
4&1&15.6344&16.0663 & 10&1&37.8824&38.0495 & 13&8&58.7789&58.9170 \\
4&2&17.0296&17.4771 & 10&2&39.2909&39.4603 & 13&9&60.1882&60.3278 \\
4&3&18.4237&18.8880 & 10&3&40.6994&40.8712 & 13&10&61.5974&61.7386 \\
4&4&19.8165&20.2988 & 10&4&42.1078&42.2820 & 13&11&63.0066&63.1495 \\
5&0&17.9873&18.3193 & 10&5&43.5161&43.6928 & 13&12&64.4158&64.5603 \\
5&1&19.3889&19.7302 & 10&6&44.9244&45.1037 & 14&0&51.1764&51.2941 \\
5&2&20.7900&21.1410 & 10&7&46.3326&46.5145 & 14&1&52.5861&52.7049 \\
5&3&22.1906&22.5518 & 10&8&47.7407&47.9254 & 14&2&53.9958&54.1158 \\
5&4&23.5906&23.9627 & 10&9&49.1487&49.3362 & 14&3&55.4055&55.5266 \\
6&0&21.7073&21.9832 & 10&10&50.5567&50.7471 & 14&4&56.8151&56.9375 \\
6&1&23.1118&23.3940 & 11&0&40.1527&40.3025 & 14&5&58.2248&58.3483 \\
6&2&24.5160&24.8049 & 11&1&41.5617&41.7133 & 14&6&59.6343&59.7591 \\
6&3&25.9199&26.2157 & 11&2&42.9706&43.1242 & 14&7&61.0439&61.1700 \\
6&4&27.3235&27.6265 & 11&3&44.3795&44.5350 & 14&8&62.4534&62.5808 \\
6&5&28.7267&29.0374 & 11&4&45.7884&45.9459 & 14&9&63.8630&63.9917 \\
6&6&30.1295&30.4482 & 11&5&47.1972&47.3567 & 14&10&65.2724&65.4025 \\
7&0&25.4109&25.6470 & 11&6&48.6059&48.7676 & 14&11&66.6819&66.8133 \\
7&1&26.8172&27.0579 & 11&7&50.0146&50.1784 & 14&12&68.0913&68.2242 \\
7&2&28.2232&28.4687 & 11&8&51.4232&51.5892 & 14&13&69.5007&69.6350 \\
7&3&29.6291&29.8796 & 11&9&52.8318&53.0001 & 14&14&70.9101&71.0459 \\
7&4&31.0347&31.2904 & 11&10&54.2403&54.4109 & 15&0&54.8482&54.9580 \\
7&5&32.4401&32.7013 & 12&0&43.8291&43.9664 & 15&1&56.2580&56.3688 \\
7&6&33.8453&34.1121 & 12&1&45.2383&45.3772 & 15&2&57.6679&57.7796 \\
8&0&29.1045&29.3109 & 12&2&46.6476&46.7881 & 15&3&59.0777&59.1905 \\
8&1&30.5119&30.7218 & 12&3&48.0568&48.1989 & 15&4&60.4875&60.6013 \\
8&2&31.9191&32.1326 & 12&4&49.4660&49.6097 & 15&5&61.8973&62.0122 \\
8&3&33.3261&33.5434 & 12&5&50.8752&51.0206 & 15&6&63.3070&63.4230 \\
8&4&34.7331&34.9543 & 12&6&52.2843&52.4314 & 15&7&64.7168&64.8338 \\
8&5&36.1399&36.3651 & 12&7&53.6933&53.8423 & 15&8&66.1265&66.2447 \\
8&6&37.5465&37.7760 & 12&8&55.1023&55.2531 & 15&9&67.5362&67.6555 \\
8&7&38.9530&39.1868 & 12&9&56.5113&56.6639 & 15&10&68.9458&69.0664 \\
8&8&40.3593&40.5976 & 12&10&57.9202&58.0748 & 15&11&70.3555&70.4772 \\
9&0&32.7915&32.9748 & 12&11&59.3291&59.4856 & 15&12&71.7651&71.8881 \\
9&1&34.1995&34.3856 & 12&12&60.7379&60.8965 & 15&13&73.1747&73.2989 \\
9&2&35.6075&35.7965 & 13&0&47.5035&47.6302 & 15&14&74.5843&74.7097 \\
\hline
\end{tabular}
\caption{Volumes $\Vol(g,k)$ and limiting volumes $\Vol^{*}(g,k)$ (for the definition see Section~\ref{sec:asymptotics}) for all admissible pairs $(g,k)$ with $2\leqslant g\leqslant 15$.}
\label{tab:volumes}
\end{table}

\end{document}